\documentclass[a4paper]{article}

\usepackage{amsmath,amssymb,amsfonts,amsthm}
\usepackage{xspace}

\ifx\xspace\undefined
  \PackageWarning{burdges}{You should be using the xspace package}
\fi

\def\newmathbf#1{\expandafter\newcommand\csname #1\endcsname{\mathbf{#1}}}
\def\newmathit#1{\expandafter\newcommand\csname #1\endcsname{\mathit{#1}}}
\def\newmathrm#1{\expandafter\newcommand\csname #1\endcsname{\mathrm{#1}}}
\def\newmathsf#1{\expandafter\newcommand\csname #1\endcsname{\mathsf{#1}}}
\def\newmathtt#1{\expandafter\newcommand\csname #1\endcsname{\mathtt{#1}}}

\def\newmathbb#1{\expandafter\newcommand\csname #1\endcsname{\mathbb{#1}}}
\def\newmathcal#1{\expandafter\newcommand\csname #1\endcsname{\mathcal{#1}}}
\def\newmathscr#1{\expandafter\newcommand\csname #1\endcsname{\mathscr{#1}}}
\def\newmathfrak#1{\expandafter\newcommand\csname #1\endcsname{\mathfrak{#1}}}

\def\newmathop#1{\expandafter\newcommand\csname #1\endcsname{\operatorname{#1}}}
\def\newmathoplim#1{\expandafter\newcommand\csname #1\endcsname{\operatorname*{#1}}}

\def\newmathbin#1{\expandafter\newcommand\csname #1\endcsname{\mathbin{\mathrm{#1}}}}
\def\newmathrel#1{\expandafter\newcommand\csname #1\endcsname{\mathrel{\mathrm{#1}}}}

\newcommand*\field[1]{\mathbb{#1}} % propser uses \R

\newcommand*\normal{\lhd}

\newcommand{\notvsim}{\not\kern-0.07cm\mid\!\sim}

\def\con#1=#2(#3){#1\equiv#2\pmod{#3}}

\def\myo{\ifmmode\circ\else$^\circ$\xspace\fi}  \let\o=\myo
\def\opri{\ifmmode^{\circ\prime}\else$^{\circ\prime}$\xspace\fi}

\newcommand*\comp{\mathbin{\raise.29ex\hbox{$\scriptscriptstyle\circ$}}}
\newcommand*\contr{\mathbin{\raise.29ex\hbox{$\lrcorner$}}}
\newcommand*\bigdot{{\mathbin{\raise.29ex\hbox{$\scriptscriptstyle{\bullet}$}}}}

\newcommand*\ov[1]{\overline{#1}}

\newcommand*\Setst[2]%
        {\left\{\,#1\vphantom{#2} \;\right|\left. #2 \vphantom{#1}\,\right\}}
\newcommand*\gen[1]{\langle #1 \rangle}
\newcommand*\Genst[2]%
        {\left<\,#1\vphantom{#2} \;\right|\left. #2 \vphantom{#1}\,\right>}

\newcommand*\Z{{\field Z}}

\newcommand*\F{{\field F}}

\newcommand*\bP{{\mathbb P}}

\newmathop{Sym}
\newmathop{Alt}

\newmathop{GL}
\newmathop{SL}
\newmathop{PSL}
\newmathop{Sp}
\newmathop{PSp}
\newmathop{SO}
\newmathop{PSO}
\newmathop{G}
\newmathop{E}
\newmathop{id}
\newmathop{Id}
\newmathop{diag}
\newmathop{sgn}
\newmathop{lcm}
\newmathop{tr}
\newmathop{Tr}
\newmathop{rk}
\newmathop{rank}
\newmathop{mult}
\newmathop{pr}
\newmathop{n}
\newmathop{m}

\newmathop{ad}
\newmathop{Ad}
\newmathop{Der}
\newmathop{Map}
\newmathop{Hom}
\newmathop{coker}

\newmathop{Aut}
\newmathop{Inn}
\newmathop{Gal}
\newmathop{Pic}
\newmathop{Stab}
\newmathop{Supp}

\newmathop{Ann}
\newmathop{Ass}
\newmathop{End}
\newmathop{Ext}
\newmathop{Tor}
\newmathop{Spec}
\newmathop{Proj}

\newmathop{Th}
\newmathop{dcl}
\newmathop{acl}
\newmathop{ord}
\newmathop{tp}
\newmathop{stp}
\newmathop{qftp}
\newmathrm{eq} % used as an exponent

\newcommand\Altinel{Alt\i nel\xspace}

\newcommand\Prufer{Pr\"ufer\xspace}

\providecommand\empty{}
\newcommand\twophrases[3]{\def\tmparg{#3}\ifx\empty\tmparg{}#1\else{}#2#3\fi}

\providecommand\Toappear[1][]{\twophrases{To appear}{To appear in }{#1}}

\providecommand\Preprint[1][]{Preprint available}

\def\secref#1{\S\ref{#1}}
\def\clref#1{Claim~\ref{#1}}

\ifx\hyperbaseurl\undefined\relax\else
  \hyperbaseurl{http://www.math.rutgers.edu/~burdges/}
\fi

\newcommand*{\hrefwww}[1]{\href{http://www.#1}{\nolinkurl{www.#1}}}

\theoremstyle{plain} %% requires amsthm

\ifx\thesection\undefined
  \newtheorem{lemma}{Lemma}
\else
  \ifx\thelemma\undefined
    \newtheorem{lemma}{Lemma}[section]
  \else
  \fi
\fi

\newtheorem{fact}[lemma]{Fact}

\newtheorem*{rawnamedtheorem}{\therawnamedtheorem}
\newcommand{\therawnamedtheorem}{\error}
\newenvironment{namedtheorem}[1]{\renewcommand{\therawnamedtheorem}{#1}
   \begin{rawnamedtheorem}}
  {\end{rawnamedtheorem}}
\newtheorem{Claim}[lemma]{Claim}

\theoremstyle{definition} %% requires amsthm

\newtheorem{hypothesis}[lemma]{Hypothesis}

\newtheorem*{rawnameddefinition}{\therawnameddefinition}
\newcommand{\therawnameddefinition}{\error}

\theoremstyle{remark} %% requires amsthm

\makeatletter
\thm@notefont{\fontseries\mddefault\upshape}%
\makeatother

   {\begin{proof}[#1]}{\end{proof}}
\newenvironment{verification}[1][Verification]%
   {\begin{proof}[#1]}{\end{proof}}

\newcounter{hypothesis}
   {\begin{list}{\alph{hypothesis}.}{\usecounter{hypothesis}}}{\end{list}}
\newcounter{conclusion}
\newenvironment{conclusions}%
   {\begin{list}{\arabic{conclusion}.}{\usecounter{conclusion}}}{\end{list}}
   {\begin{conclusions}}{\end{conclusions}}
\newcommand\mathperiod{.}

\def\-#1{\overline{#1}}

\def\leqs{\leqslant}
\def\geqs{\geqslant}

\title{Uniqueness Cases in Odd Type Groups of Finite Morley Rank}
\author{Alexandre V.\ Borovik%
\thanks{The first author completed his work on the paper during his visit
to Institut Giscard Desargues, Universit\'{e} Lyon 1, in April 2003.}\\
School of Mathematics, The University of Manchester\\
PO Box 88, Sackville St., Manchester M60 1QD, England\\
{\tt alexandre.borovik@umist.ac.uk}
\and
Jeffrey Burdges\thanks{Partially supported by DFG grant Te 242/3-1.}\\
Mathematisches Institut, Universit\"at W\"urzburg\\
Am Hubland, D-97074 W\"urzburg, Germany\\
{\tt burdges@math.rutgers.edu}
\and
Ali Nesin\thanks{Partially supported by the London Mathematical Society grant 4523.}\\
Mathematics Department,  Istanbul Bilgi University\\
Ku\c{s}tepe \c{S}i\c{s}li, Istanbul, Turkey\\
{\tt anesin@bilgi.edu.tr}
}
\begin{document}

\maketitle

\begin{abstract}
There is a longstanding conjecture, due to Gregory Cherlin and Boris Zilber,
that all simple groups of finite Morley rank are simple algebraic groups.
One of the major theorems in the area is Borovik's trichotomy theorem.
The ``trichotomy'' here is a case division of the generic minimal
counterexamples within {\it odd type}, i.e.~groups whose Sylow$^\circ$
2-subgroup is large and divisible.  The so-called uniqueness case the
trichotomy theorem is the existence of a proper 2-generated core.
It is our goal to drive presence of a proper 2-generated core to a
contradiction; and hence bound the complexity of the Sylow$^\circ$
2-subgroup of a minimal counterexample to the Cherlin-Zilber conjecture.
This paper shows that the group in the question is a minimal connected
simple group and has a strongly embedded subgroup, a far stronger
uniqueness case.  As a corollary, a tame counterexample to the
Cherlin-Zilber conjecture has \Prufer rank at most two.

% There is a longstanding conjecture, due to Gregory Cherlin and Boris Zilber,
% that all simple groups of finite Morley rank are simple algebraic groups.
% One of the major remaining obstacles to this conjecture is odd type
% $K^*$-groups with high \Prufer 2-rank and a proper 2-generated core.
% This paper shows that the group in the question is a minimal connected
% simple group and has a strongly embedded subgroup.

\footnotetext{Mathematics Classification Numbers 03C60, 20G99}
\end{abstract}

\section{Introduction}\label{sec:Intro}

This paper relates to the algebraicity conjecture for simple groups of
finite Morley rank, also known as the Cherlin-Zilber conjecture, which
states that all simple groups of finite Morley rank are simple algebraic
groups over an algebraically closed field.  As with most of the recent
work on this conjecture, the present article seeks to transfer ideas
from the classification of finite simple groups.
% to groups of finite Morley rank.

It is now common practice to divide the Cherlin-Zilber conjecture into
different cases depending on the nature of the connected component of the
Sylow 2-subgroup, or Sylow\o 2-subgroups (cf. \S\ref{sec:Sylow}).
We shall be working with groups whose Sylow\o 2-subgroup is divisible
and non-trivial, or {\em odd type} groups.  Prior to \cite{BuPhd},
the main theorems in the area of odd type groups are Borovik's Trichotomy 
Theorem \cite{Bo95} and the Generic Identification Theorem \cite{BB01}.
Together, these two results prove the following.

\begin{namedtheorem}{Tame Trichotomy theorem}
Let $G$ be a simple {\em tame} $K^*$-group of finite Morley rank and
odd type.  Then $G$ is either a Chevalley group over an algebraically
closed field of characteristic not 2, or has normal 2-rank $\leqs 2$,
or has a proper 2-generated core.
\end{namedtheorem}  % newtrichotomy is often Bo03

Here a group is said to be tame if it does not involve a field of finite Morley
rank with a proper infinite definable subgroup of it's multiplicative group.
Such fields are presently believed to exist in characteristic zero
\cite{Po01b}.  Hence the tameness assumption must eventually be
removed.

In this paper, we analyze groups with proper 2-generated cores
(see \S\ref{sec:StrEmbDef} for the definition), and
drive them towards exceptional minimal connected simple configurations
which should eventually turn out to be contradictory.
In \cite{CJ01}, Cherlin and Jaligot show that the \Prufer 2-rank of
a tame minimal connected simple group is at most 2.  In light of this result,
and the Tame Trichotomy, the present paper shows the following.

\begin{namedtheorem}{Tame Generic Case}
A tame minimal counterexample to the algebraicity conjecture has
\Prufer 2-rank at most 2.
\end{namedtheorem}

It is our near term goal to eliminate the need for tameness in the above theorem.
In \cite{BuPhd}, tameness is removed from the tame trichotomy above,
and the present paper will make no use of tameness either, so all important
applications of tameness now lie within \cite{CJ01}.  For this reason,
our results below will push beyond establishing that the group is minimal
connected simple, and attempt to provide tools for the analysis of minimal
connected simple groups, without tameness.  In particular, we will show that
the Sylow 2-subgroup is connected, and that $G$ has a strongly embedded
subgroup.  Our results are summarized as follows.

\begin{namedtheorem}{Strong Embedding Theorem}
Let $G$ be a simple $K^*$-group of finite Morley rank and odd type with
normal 2-rank $\geqs 3$ and \Prufer 2-rank $\geq2$.  Let $S$ be a Sylow
2-subgroup of $G$.  Suppose that $G$ has a proper 2-generated core
$M = \Gamma_{S,2}(G) < G$.  Then the following hold.
\begin{enumerate}
\item $G$ is a minimal connected simple group, i.e.\ 
  all proper definable connected subgroups are solvable.
\item $M$ is strongly embedded.
\item $B := M^\o$ is a Borel subgroup
\item $S$ is connected.
%referee's request
\item[5]  $N_G(B) = M$.
\item[6]  $I(B \cap B^g) = \emptyset$ for any $g\notin M$.
\item[7]  $\bigcup B^G$ is generic in $G$.
\end{enumerate}
\end{namedtheorem}

Burdges, Cherlin, and Jaligot will eliminate this configuration in \cite{BCJ},
thus replicating the main result of \cite{CJ01}.

The notions of both 2-generated core and strongly embedded subgroup
arise as so-called {\em uniqueness cases} in finite group theory.
These subgroups both exhibit a black hole property reminiscent of a normal
subgroup; and they seem similar when we compare Fact~\ref{strembdef}
below with Lemma~\ref{stremb_CA} or \clref{stremb_CM} of \secref{sec:StrEmb}.
Strong embedding, however, is far more powerful and has global consequences
(see Fact \ref{stremb_invconj}).
Our proof of the fact that $G$ is a minimal connected simple group will
involve passing through strong embedding to obtain a contradiction under
the assumption that $B$ is non-solvable.

In bridging the gap between 2-generated cores and strong embedding,
we employ the theory of Carter subgroups and make use of a result due
to Olivier Fr{\'e}con (Fact~\ref{Carter_Sylow}) in the final stage of
the argument.

\section{Background}\label{sec:BG}

We now recall essential facts about groups of finite Morley rank.
The standard reference for our basic facts is \cite{BN}.  Some of that
material will be used without explicit mention.

A group of finite Morley rank is {\em connected} if it contains no
proper definable subgroup of finite index.  We will refer to maximal
connected solvable subgroups of a group of finite Morley rank as
{\em Borel} subgroups.

We define the {\em 2-rank} $m_2(G)$ of a group $G$ to be the maximum
rank of its elementary abelian 2-subgroups.  Also, the {\em \Prufer 2-rank}
$\pr_2(G)$ is the maximum rank of its \Prufer 2-subgroups $\Z(2^\infty)^k$,
and the {\em normal 2-rank} $n_2(G)$ is the maximum rank of a normal
elementary abelian 2-subgroup of Sylow 2-subgroup of $G$.  These ranks
must all be finite for subgroups of an odd type group of finite Morley rank.

We define the {\em odd part} $O(G)$ of a group $G$ of finite Morley rank
to be the maximal definable connected normal $2^\perp$-subgroup of $G$.
The subgroup $O(G)$ is well-defined by the following exercise from \cite{BN}.

\begin{fact}[Exercise 11 on page 93 of \cite{BN}]\label{nodeftorsion}
Let $G$ be a group of finite Morley rank and let $H \normal G$ be a
definable subgroup.  Let $x\in G$ be an element such that $\bar{x}\in G/N$
is a $p$-element.  Then $x H$ contains a $p$-element.
\end{fact}

\subsection{Sylow and Carter subgroups}\label{sec:Sylow}

We provide a basic notion of ``characteristic'' for groups of finite
Morley rank as follows.

Let $S$ be a Sylow 2-subgroup of a group $G$ of finite Morley rank.
By \cite{BP} (see also Lemma 10.8 of \cite{BN}),
 $S^\o = B * T$ is a central product of a definable connected
nilpotent subgroup $B$ of bounded exponent and of a 2-torus $T$,
i.e.~$T$ is a divisible abelian 2-group.

The group $G$ is said to have {\em odd type} if $B = 1$ and $T \neq 1$.
This notion is well-defined because the Sylow 2-subgroups of a
group of finite Morley rank are conjugate by \cite{BP,PW93}
(see also Theorem 10.11 of \cite{BN})
The following two corollaries of conjugacy, known as a ``Frattini
argument'' and a ``fusion control lemma'' respectively, will be useful.

\begin{fact}[Corollary 10.12 of \cite{BN}]\label{Sylow_frattini}
Let $G$ be a group of finite Morley rank, let $N \normal G$ be a
definable subgroup, and let $S$ be a characteristic subgroup of the
Sylow 2-subgroup of $N$.  Then $G = N_G(S) N$.
\end{fact}

\begin{fact}[\S10.6.1 of \cite{BN}]\label{Sylow_fusion}
Let $G$ be a group of finite Morley rank and odd type.
Let $S$ be a Sylow 2-subgroup of $G$.
Then $N_G(S^\o)$ controls fusion in $C_S(S^\o)$, i.e.~two elements of
$C_S(S^\o)$ which are $G$-conjugate are in fact $N_G(S^\o)$-conjugate.
\end{fact}

A useful property of Sylow 2-subgroups is that they can be lifted:

\begin{fact}[\cite{PW00}; Corollary 1.5.5 of \cite{Wag}]
\label{Sylow_lift}
Let $G$ be a group of finite Morley rank and let $N$ be a normal subgroup
of $G$.  Then the Sylow 2-subgroups of $G/N$ are the images of the Sylow
2-subgroups of $G$.
\end{fact}

\begin{fact}[Theorem 9.29 of \cite{BN}; see also Corollary 7.15 of \cite{Fre00a}]
\label{Sylow_con}
Let $G$ be a connected solvable group of finite Morley rank.
Then the Sylow $p$-subgroups of $G$ are connected.
\end{fact}

% \subsection{Carter subgroups}\label{sec:Carter}

Let $G$ be a group of finite Morley rank.
A definable subgroup $C \leqs G$ which is nilpotent and self-normalizing in $G$
 is called a {\em Carter subgroup} of $G$.

The following result is a summary, in order, of \cite[Proposition 3.2]{Fre00a},
\cite[Corollary 4.8]{Fre00a}, \cite[Theorem 5.5.12]{Wag,Wag94},
and \cite[Corollary 7.15]{Fre00a}.

\begin{fact}
\label{Carter_cleanup}\label{Carter_con}\label{Carter_exists}
\label{Carter_Sylow}\label{Carter_conj}\label{Carter_frattini}
Let $H$ be a connected solvable group of finite Morley rank. Then the following hold.
\begin{itemize}
\item[(1)] $H$ has a Carter subgroup.
\item[(2)] The Carter subgroups of $H$ are the definable nilpotent
subgroup of $H$ with $N_H^\o(C) = C$.
In particular, Carter subgroups of $H$ are connected.
\item[(3)] The Carter subgroups of $H$ are $H$-conjugate.
As a corollary, we have the Frattini argument: if $H$ is a definable connected
normal subgroup of a group $G$ of finite Morley rank, and $C$ is a
Carter subgroup of $H$, then $G = N_G(C) H$.
\item[(4)] Let $R$ be a Sylow $p$-subgroup of $H$.
Then $N_H(R)$ contains a Carter subgroup of $H$.
\end{itemize}
\end{fact}

% \subsection{Genericity}\label{sec:Gen}
%
% We also have a few lemmas which Eric Jaligot and Gregory Cherlin
% found useful in their analysis of tame minimal connected simple
% groups; these elaborate on a genericity argument used to study bad
% groups \cite[Ch. 13]{BN}.

\subsection{Algebraic groups and $K$-groups}\label{sec:Kgrp}

A group $G$ will be called {\em quasi-simple} if $G=G'$ and $G/Z(G)$ is simple.
The group $G$ will be called {\em semi-simple} if $G=G'$ and $G/Z(G)$ is
{\em completely reducible}, i.e.~$G/Z(G)$ is a direct sum of finitely many simple
subgroups.   So quasi-simple groups are semi-simple.

% We will not follow the tradition of using the phrase ``simple algebraic
% group'' to refer to a connected algebraic group $G$ with $Z(G)$ finite
% and $G/Z(G)$ simple as an abstract group.  Instead, we will use the
% phrase quasi-simple algebraic group and depend on Fact~\ref{centralext}
% below to keep the center finite.

We will need the following results from the classification of quasi-simple
algebraic groups.

\begin{fact}\label{alggrpmin} % [\cite{???}]
The only quasi-simple algebraic groups over an algebraically closed field
$F$ without proper definable quasi-simple subgroups are $\SL_2(F)$ and $\PSL_2(F)$.
\end{fact}

\begin{fact}[Theorem 8.4 of \cite{BN}] \label{autalg}
Let $G \rtimes H$ be a group of finite Morley rank where $G$ and
$H$ are definable, $G$ an infinite quasi-simple algebraic group over an
algebraically closed field, and $C_H(G)$ is trivial.  Then, viewing $H$
as a subgroup of $\Aut(G)$, we have $H\leqs \Inn(G)\Gamma$, where
$\Inn(G)$ is the group of inner automorphisms of $G$ and $\Gamma$
is the group of graph automorphisms of $G$, relative to a fixed choice
of Borel subgroup $B$ and maximal torus $T$ contained in $B$.
\end{fact}

A group $G$ of finite Morley rank is called a {\em $K$-group} if every
connected definable simple section of $G$ is a Chevalley group over an
algebraically closed field.  We shall also call a group $G$ of finite
Morley rank a {\em $K^*$-group} if every proper definable section is
a $K$-group.  Clearly, a {\em minimal} non-algebraic connected simple
group of finite Morley rank will be a $K^*$-group.  We also observe
that $O(H)$ is solvable if $H$ is a $K$-group, since simple algebraic
groups contain involutions.

A quasi-simple subnormal subgroup of a group $G$ is referred to as
a {\em component} of $G$.

\begin{fact}[\cite{Bel87,Ne91}; see also \S7.4 of \cite{BN}]\label{components}
Let $G$ be a group of finite Morley rank.  Then the components of $G$
are definable subgroups, and there are only finitely many of them.
Furthermore, $G$ acts by conjugation on the set of components
(see Lemma 7.12ii of \cite{BN}).
\end{fact}

The subgroup $L(G)$ generated by the components of $G$ is now
definable, being the setwise product of the components.  We will
refer to $L(G)$ as the {\em layer} of $G$ and define $E(G) = L^\o(G)$.

\begin{fact}[\cite{AC99}] \label{centralext}
A group of finite Morley rank which is a perfect central extension of
a quasi-simple algebraic group over an algebraically closed field is
an algebraic group and has finite center.
\end{fact}

We define the {\em Fitting subgroup} $F(G)$ of $G$, to be the subgroup
generated by all the normal nilpotent subgroups of $G$.  The Fitting subgroup is
nilpotent and definable \cite{Bel87,Ne91} (see also \cite[Theorem 7.3]{BN}).

\begin{fact} \label{511} % \label{Bo95}
Let $G$ be a connected $K$-group of odd type. Then $G/O(G)$ is
isomorphic to a central product of quasi-simple algebraic groups
over algebraically closed fields of characteristic not 2 and of
a definable connected abelian group. In particular, if $\-G=G/O(G)$
then $\-G=F(\-G)E(\-G)$ and $F(\-G)$ is an abelian group.
\end{fact}

\begin{proof}
The ``in particular'' part of the statement is \cite[Theorem~5.9]{Bo95}.
By  definition, $E(G) = L_1 * \cdots * L_k$ is a central product of
connected quasi-simple groups.  Since $G$ is a $K$-group, each $L_i$ is
a perfect central extension of a Chevalley group over an algebraically
closed field.  Now the result follows from Fact~\ref{centralext}.
\end{proof}

A {\em Klein four-group}, or just {\em four-group} for short, is a
group isomorphic to $\Z/2\Z \oplus \Z/2\Z$.  We will use the notation
$H^\# = H \setminus \{1\}$ to denote the set of non-identity elements
of a group $H$.

The following generation principle for $K$-groups will be used frequently.

\begin{fact}[Theorem 5.14 of \cite{Bo95}]\label{Kgrp_gen}
Let $G$ be a connected $K$-group of finite Morley rank and odd type.
Let $V$ be a four-subgroup acting definably on $G$.  Then
$$ G = \gen{ C^\o_G(v) \mid v\in V^\# } $$
\end{fact}

\section{Uniqueness subgroups}\label{sec:StrEmbDef}

We first discuss the notions of 2-generated core and strongly embedded subgroup.

A proper definable subgroup $M$ of a group $G$ of finite Morley rank
is said to be {\em strongly embedded} if $I(M) \neq \emptyset$ and
$I(M \cap M^g) = \emptyset$ for any $g\in G \backslash M$.
Here $I(H)$ to denotes the set of involutions of $H$.
We will apply the usual criteria for strong embedding:

\begin{fact}[Theorem 9.2.1 of \cite{Gor}; see also Theorem 10.20 of \cite{BN}]%
\label{strembdef} % Fact 3.1 of \cite{Al94}
Let $G$ be a group of finite Morley rank with a proper definable
subgroup $M$.  Then the following are equivalent:
\begin{enumerate}
\item $M$ is a strongly embedded subgroup.
\item $I(M) \neq \emptyset$, $C_G(i) \leqs M$ for every $i \in I(S)$,
      and $N_G(S) \leqs M$ for some Sylow 2-subgroup $S$ of $M$.
\item $I(M) \neq \emptyset$ and $N_G(S) \leqs M$ for every non-trivial
      2-subgroup $S$ of $M$.
\end{enumerate}
\end{fact}

The following is one of the major applications of strong embedding.

\begin{fact}[Theorem 10.19 of \cite{BN}; see alse Fact 3.3 of \cite{Al96}]%
\label{stremb_invconj}
Let $G$ be a group of finite Morley rank with a proper definable strongly
embedded subgroup $M$.  Then
\begin{enumerate}
\item A Sylow 2-subgroup of $M$ is a Sylow 2-subgroup of $G$,
\item $G$ and $M$ each have only one conjugacy class of involutions,
\end{enumerate}
\end{fact}

Let $G$ be a group of finite Morley rank and let $S$ be a Sylow 2-subgroup
of $G$.  We define the {\em 2-generated core} $\Gamma_{S,2}(G)$ of $G$
to be the definable hull of the group generated by all normalizers $N_G(U)$
of all elementary abelian 2-subgroups $U \leqs S$ with $m_2(U) \geqs 2$.
As it is this last rank condition to which the ``2-generated'' is referring,
a strongly embedded subgroup would be a proper 1-generated core
by Fact \ref{strembdef}.

A priori, merely possessing a proper 2-generated core need not entail the
global consequences of Fact \ref{stremb_invconj}.  However, the following
easy consequence of Fact~\ref{Kgrp_gen} indicates that 2-generated cores
are not far from being strongly embedded.

\begin{lemma}\label{stremb_CA}
Let $G$ be a simple $K^*$-group of finite Morley rank and odd type.
Let $S$ be a Sylow 2-subgroup of $G$ and let $M = \Gamma_{S,2}(G)$
be the 2-generated core associated with $S$.
Let $A$ be an elementary abelian 2-subgroup of $M$ with $m_2(A) \geqs 3$.
Then $C^\o_G(a) \leqs M$ for any $a\in A^\#$.
\end{lemma}

\begin{proof}
Let $K= C^\o_G(a)$. Let $A_1$ be a four-subgroup of $A$ disjoint from
$\gen{a}$. Consider the $K$-group $K$ of odd type, which contains $A_1$.
By Fact~\ref{Kgrp_gen}, $K = \gen{C_K^\o(x)\mid x\in A_1^\#}$.
Now $C_K^\o(x) \leqs C_G(a,x)$ and $\gen{a, x}$ is a four-subgroup
of $S$.  Thus $K \leqs H$.
\end{proof}

\noindent This shows that 2-generated cores exhibit a kind of
``black hole'' principle, limiting communication between elements of
the subgroup $\Gamma_{S,2}(G)$ and its exterior.
% Reviewer:  What does this mean?

\begin{lemma}\label{twogencore_max_pr}
Let $G$ be a simple $K^*$-group of finite Morley rank and odd type.
Let $S$ be a Sylow 2-subgroup of $G$ and let $M = \Gamma_{S,2}(G)$
be the 2-generated core associated with $S$.  If $\pr_2(S) \geqs 3$
and $M < G$ then $B := M^\o$ is a maximal proper connected subgroup
of $G$.
\end{lemma}

\begin{proof}
Let $K < G $ be a connected group containing $B$.
Since $\pr_2(S) \geqs 3$, Fact~\ref{Kgrp_gen} and Lemma \ref{stremb_CA} yield
\[ K \leqs \langle C^\o_K(i) : i\in \Omega_1(S^\o) \rangle \leqs M
  \qedhere \]
\end{proof}

\section{Component Analysis}\label{sec:Component}

Our next few lemmas are directed toward the proof that $B$ is solvable.
The first of these will allow us to prove that $M$ is strongly embedded
when $B$ is non-solvable.

\begin{lemma}\label{stremb_nonsolv}
Let $G$ be a $K$-group of finite Morley rank and odd type with non-solvable
connected component and $i$ an involution in $G$.
Then the Sylow\o 2-subgroups of $C^\o_G(i)$ are non-trivial.
\end{lemma}

We first recall the following lemma.

\begin{fact}[Fact 3.2 of \cite{Bu03}; Fact 3.12 of \cite{BuPhd}]
\label{Cquotient}
Let $G = H \rtimes T$ be a group of finite Morley rank
 with $H$ and $T$ definable.
Suppose $T$ is a solvable $\pi$-group of bounded exponent and
$Q \normal H$ is a definable solvable $T$-invariant $\pi^\perp$-subgroup.
Then $$ C_H(T)Q/Q = C_{H/Q}(T)\mathperiod $$
\end{fact}

\begin{proof}[Proof of Lemma \ref{stremb_nonsolv}]
% By Fact \ref{Cquotient}, $C_{G^\o}(i)O(G)/O(G) = C_{G^\o/O(G)}(i)$.
% So we may assume $O(G)=1$ by Fact \ref{Sylow_lift}.
We claim that it is enough to prove the statement for $\-G=G/O(G)$.
Let $i\in I(G)$ and let us assume that we know the result for $\-G$
 and the involution $\-i$ of $\-G$.
Let $\-S$ be a nontrivial Sylow\o 2-subgroup of $C_{\-G}(\-i)$.
Since $C_G(i)/C_{O(G)}(i) \cong C_{\-G}(\-i)$ by Fact \ref{Cquotient},
 there is a nontrivial Sylow\o 2-subgroup $S$ of $C_G(i)$
 by Fact~\ref{Sylow_lift}.
Hence we can assume that $O(G^\o)=O(G)=1$. % (Fact~\ref{osolv}).

Let $i\in I(G)$. By Fact~\ref{511}, $G^\o$ is the central product of finitely
many quasi-simple algebraic groups and of a definable connected abelian
group $F:=F(G)^\o \normal G$, say $G^\o = G_1 * \cdots * G_n * F$.
Let $L= G_1*\cdots * G_n$. Since $L \neq 1$ and $i$ normalizes $L$ by
Fact~\ref{components}, we can assume that $G= L\rtimes \gen{i}$.
If $i$ swaps two of the quasi-simple components $G_j$ and $G_k$, then
$\gen{ss^i \mid s\in S}$, where $S$ is a Sylow 2-subgroup of $G_j$, is
an infinite 2-subgroup of $C_G(i)$ and we are done.  Therefore we may
assume that $i$ normalizes each component. This allows us to assume that
$L$ is just one component, i.e.\ $G=L\rtimes \gen{i}$ and $L$ is
quasi-simple algebraic.

By Fact~\ref{autalg}, we have two cases: $i$ acts on $L$ either as
an inner automorphism, or as an inner automorphism composed with a
graph automorphism, and hence $G$ is algebraic.
Since $G$ has odd type, $i$ is semisimple in $G$.
So $C^\o_G(i)$ is nontrivial and reductive by Theorem 8.1 of \cite{St2},
 and hence has an infinite Sylow 2-subgroup.
Alternatively, scrutinizing the table of centralizers of involutive
automorphisms of algebraic groups \cite[Table~4.3.1]{GLS3}
shows that they always have infinite Sylow 2-subgroups.
% Greg pointed out that we could just prove what we need here without
% all that much work, but I actually kinda like the reference to Lyons
\end{proof}

\noindent The next lemma will be used to contradict strong embedding
under the assumption that $B$ is non-solvable.

\begin{lemma} \label{nonconjugacy}
Let $G$ be a $K$-group of finite Morley rank and odd type with
non-solvable $G^\o$ and $\pr_2(G) \geqs 3$.
Let $S$ be a Sylow $2$-subgroup of $G$. Then not all the
involutions of $S^\o$ are $G$-conjugate.
\end{lemma}

Notice that the assumption $\pr_2(G) \geqs 3$ cannot be
weakened: if $K$ is an algebraically closed field of
characteristic distinct from $2$ then the group $G={\rm PSL}_3(K)$
has \Prufer 2-rank 2, and only one conjugacy class of involutions.

\medskip

\begin{proof}
Suppose toward a contradiction that the involutions of $S^\o$ are
all $G$-conjugate.  Passing to a quotient, we may suppose $O(G)=1$.

By Fact~\ref{511}, $G^\o$ is a central product of finitely many
quasi-simple algebraic groups and of a definable connected abelian
group $F$, say $G^\o = G_1 * \cdots * G_n * F$.
Let $L= G_1*\cdots * G_n$.  Since $G_1$ has an involution and
$L\normal G$, all the involutions of $G$ are in $L$.

{\bf Case 1. $Z(L)$ has an involution.}
Then all the involutions of $G^\o$ are in $Z(L)$.
Thus each $G_i$ is a quasi-simple algebraic group whose involutions
are in $Z(G_i)$.  From the classification of quasi-simple algebraic
groups (e.g.\ \cite{seitz}), it follows that $G_i \simeq \SL_2(K_i)$
for some algebraically closed $K_i$ of characteristic not 2
(see Theorem 1.12.5d of \cite{GLS3}).
Thus $L$ is a central quotient of
  $\SL_2(K_1) \times \cdots \times \SL_2(K_n)$.
Any nontrivial central quotient of
  $\SL_2(K_1) \times \cdots \times \SL_2(K_n)$
will introduce new noncentral involutions since the involution of
$Z(\SL_2(K_i))$ has a noncentral square root.
So $G^\o = \SL_2(K_1) \times \cdots \times \SL_2(K_n)$.
Since $G$ permutes the components $G_1,\ldots,G_n$ by
Fact~\ref{components}, the associated set of involutions
$\{i_1,\ldots, i_n\}$, given by $i_j\in I(G_j)$, is $G$-invariant.
So $i_1$ can not be conjugate to $i_1 i_2$ if $i_1 \neq i_2$.
Since $\pr_2(G) \geq 3$ and $\pr_2(\SL_2(F_i)) = 1$,
there are at least three components, a contradiction.

{\bf Case 2. $Z(L)$ has no involutions.}
Passing to a quotient by Fact \ref{Sylow_lift},
we can assume without loss of generality that $Z(L)=1$ and that
each $G_i$ is an algebraic group over an algebraically closed field
of characteristic not 2 which is simple as an abstract group.
So $L = G_1 \times \cdots \times G_n$.
Then $S = S_1 \times \cdots \times S_n$
 with $S_i$ a Sylow 2-subgroup of $G_i$.
If $n \geq 2$ then an involution in $S_1$ cannot be conjugate to
a product of involutions from $S_1$ and $S_2$, so $n=1$.
%%
% As in the previous case, each involution belongs to a component, so $n=1$.
%%
%Let $S$ be a Sylow 2-subgroup of $L$.  Let $A$ be the subgroup of
%$S^\o$ generated by all involutions in $S^\o$.
%Since $Z(L) = 1$ and $L = G_1 \times \cdots \times G_n$,
%it follows that $A = (A \cap G_1) \times \cdots \times (A \cap G_n)$.
%By Fact~\ref{Sylow_fusion}, $N_G(S^\o)$ acts transitively on the
%involutions of $A$. The action of $G$ by conjugation permutes the
%set of components $\{G_1,\dots, G_n\}$ by Fact~\ref{components}.
%So if $n \geqs 2$, the action of $N_G(S^\o)$ on $A$ is
%imprimitive: it permutes the subgroups $A \cap G_1, \dots, A \cap
%G_n$. Thus if $a_1\in A^\#\cap G_1$ and $a_2\in A^\#\cap G_2$,
%then $a_1$ and $a_1a_2$ cannot be conjugate.  Hence $n=1$.
%%
Thus $G$ acts transitively on the involutions of the simple
algebraic group $L = G_1$.
Since $\pr_2(L) \geq 3$, there are two involutions $t,s \in L$
with $C^\o(s) \not\cong C^\o(t)$ by Table~4.3.1 of \cite{GLS3}.
So the result follows.
\end{proof}

The following lemmas will be used to show that $G$ is a minimal
connected simple group once we have the solvability of $B := M^\o$.
The first is a lifting lemma for 2-generated cores and the second is a
structural result about a group of the form $\PSL_2(K)$.

\begin{lemma} \label{lifting}
Let $G$ be a group of finite Morley rank and odd type.
Let $S$ be a Sylow $2$-subgroup of $G$.
Let $\bar{\quad}$ denote ``image in the quotient $G/O(G)$.''
Then $$\-{\Gamma_{S,2}(G)}=\Gamma_{\-S,2}(\-G)$$
\end{lemma}

\begin{proof}
For any four-group $A\leq S$, the image $\ov A$ is still a four-group.
So the left hand side is a subgroup of the right hand side.
To prove the reverse inclusion, it is enough to show that,
for any four-subgroup $E$ of $\ov S$, we have a four-subgroup $A$ of $S$
such that $\ov A = E$ and $N_{\ov G}(\ov A) \leq \ov{N_G(A)}$.

Let $E$ be a four-subgroup of $\ov S$ and let $X$ be the full preimage
of $E$ in $G$.  Since $E \leq \ov S$, we have $X\leq SO(G)$.
Let $A$ be a Sylow 2-subgroup of $X$.  By Fact \ref{Sylow_lift},
$\ov A = E$, so $X = A O(G)$ and $A \cong E$.
Since $A\leq X \leq S O(G)$ and $S$ is a Sylow 2-subgroup of $S O(G)$,
we may assume that $A\leq S$ by conjugating by an element of $O(G)$.
Since $A$ is a Sylow 2-subgroup of $A O(G)$,
$N_G(A O(G)) \leq N_G(A) O(G)$ by Fact \ref{Sylow_frattini}.
So $N_{\ov G}(\ov A) \leq \ov{N_G(A)}$, as desired.
\end{proof}

\begin{lemma} \label{notpsl2}
The connected component of a 2-generated core of $\PSL_2(K)$,
where $K$ is an algebraically closed field of characteristic distinct
from 2, is non-solvable.
\end{lemma}

Notice that it follows from Poizat \cite{Po01a} that $\PSL_2(K)$
coincides with its 2-generated core, although we do not need the
full strength of this result.

\medskip

\begin{proof}
Let $T$ be the standard maximal torus of $G=\PSL_2(K)$
(that consists of diagonal elements modulo the center of $\SL_2(K)$).
Let $S$ be a Sylow 2-subgroup of $G=\PSL_2(K)$ such that
$S^\o \leqs T$.  Then $S=S^\o \rtimes \gen{w}$ for some
$w\in I(N_G(T) \setminus T)$. % by Fact~\ref{???}.
Since $w$ inverts $T$, $w S^\o$ consists entirely of involutions
 and $S$ is generated by its involutions.
Let $z$ be the unique involution of $Z(S)\leqs S^\o$.
Let $M = \Gamma_{S,2}(G)$.
For any involution $t\neq z$ of $S$, $t$ belongs to the
four-subgroup $\gen{z,t}$ of $S$, so $S \leqs M$.

Now recall that $G$ is the automorphism group of the projective
line $\bP^1$ over the field $K$.
Since $z$ and $t$ are involutions and the characteristic is not 2,
they have two fixed points each, which we label $z_1,z_2$ and $t_1,t_2$,
respectively.
Since $t$ commutes with $z$, they stabilize one another's fixed points.
Since $z \neq t$, we have $z_1 \neq t_1$ and $z_1 \neq t_2$.
Also $z_1^t = z_2$ and $z_2^t = z_1$.
Since $G$ acts sharply 3-transitively on $\bP^1$, % REFERENCE
there is an $r\in G$ such that
$z_1^r =  t_1$, $z_2^r = t_2$, and $t_1^r = z_1$.
Since the pointwise stabilizer of $t_1$ and $t_2$ is isomorphic to $K^*$,
there is only one involution fixing these two points, and thus $z^r = t$.
Since $t^r$ commutes with $z^r = t$, $t$ stabilizes the fixed point set of $t^r$.
Since $t^r$ fixes $z_1 = t_1^r$, and $z_1^t = z_2$,
we find that $t^r$ fixes $z_2$ too, and thus $t^r = z$.
Hence $r$ normalizes $\gen{z,t}$ and $r\in M$.
Now $t \in S^{\o r^{-1}}$ since $z\in S^\o$.
So $\gen{z,t} \leq M^\o$.

Suppose towards a contradiction that $M^\o$ is solvable.
Then $\pr_2(M^\o) \geq 2$ by Fact \ref{Sylow_con},
 a contradiction.
\end{proof}

\section{Proof of the Strong Embedding Theorem}\label{sec:StrEmb}

Let $G$ be a simple $K^*$-group of finite Morley rank and odd type with
normal 2-rank $\geqs 3$ and \Prufer 2-rank $\geq2$.
Let $S$ be a Sylow 2-subgroup of $G$.  Suppose that $G$ has
 a proper 2-generated core $M = \Gamma_{S,2}(G) < G$.
We proceed by first establishing that $G$ is a minimal connected simple group,
and then showing that $S$ is connected, which can be used to prove
strong embedding of $M$.

Let $E \normal S$ be an elementary abelian 2-subgroup with $m_2(E)\geq3$.

\begin{Claim}\label{stremb_E1}
For every $i\in I(S)$, $C_E(i)$ contains a four-group.
\end{Claim}

\begin{verification}
Since $E$ is normal in $S$, the involution $i$ induces a linear transformation
of the $\F_2$-vector space $E$.  Since $m_2(E) > 2$, the Jordan canonical form
of $i$ cannot consist of a single block, so there are at least two eigenvectors.
Since the eigenvalues associated to these eigenvectors must have order 2, the
eigenvalues must both be 1, as desired.
\end{verification}

\begin{Claim}\label{stremb_CconM}
$C^\o_G(i) \leq M$ for every $i\in I(M)$.
\end{Claim}

\begin{verification}
We may assume that $i\in I(S)$ after conjugation. % by Fact \ref{Sylow}
By \clref{stremb_E1}, there is a four-group $E_1 \leq E$ centralized
by $i$.  Thus either $E$ or $\langle E_1,i \rangle$ is an elementary
abelian 2-group of rank at least three which contains $i$.
By Lemma~\ref{stremb_CA}, $C^\o_G(i) \leq M$.
\end{verification}

\begin{Claim}\label{stremb_CM}
$C_G(i) \leq M$ for any $i\in I(M)$ for which $C^\o_M(i)$ has an infinite
Sylow 2-subgroup.
\end{Claim}

\begin{verification}
Let $R$ be a Sylow\o\ 2-subgroup of $C^\o_G(i)$.
% Tuna objects to the \o here since they are all connected
We may assume that $\langle R, i \rangle \leq S$ after conjugation.% by Fact~\ref{Sylow}.
We claim that $N_{C_G(i)}(R) \leq M$.
If $i\notin S^\o$ then $m_2(\langle \Omega_1(R),i \rangle) \geq 2$,
so  %%% (if either $\pr_2(S)=1$ or $\pr_2(S)>1$)
$$ N_{C_G(i)}(R) \leq N_G(\langle \Omega_1(R),i \rangle) \leq M $$
If $i\in S^\o$ then $m_2(\Omega_1(R))\geq 2$ since $\pr_2(G) \geq 2$,
so $$ N_{C_G(i)}(R) \leq N_G(\Omega_1(R)) \leq M $$
Now Fact~\ref{Sylow_frattini} and \clref{stremb_CconM} yield
\[ C_G(i) = C^\o_G(i) N_{C_G(i)}(R) \leq M \qedhere \]
\end{verification}  %% solvable Fratini for $p \neq 2$?

\begin{Claim}\label{stremb_Msolv}
$B := M^\o$ is solvable.
\end{Claim}\nobreak

\begin{verification}
Suppose toward a contradiction that $B$ is non-solvable.
Then, by Lemma~\ref{stremb_nonsolv}, for every involution $i \in M$,
the Sylow 2-subgroups of $C_M(i)$ are infinite.
By Claim~\ref{stremb_CM}, $C_G(i) \leqs M$.
So $M$ is strongly embedded by Fact~\ref{strembdef}, and
 any two elements of $E^\#$ are $G$-conjugate by Fact~\ref{stremb_invconj}.

We observe that $N_G(S^\o) \leqs M$ since $\pr_2(S) \geqs 2$.
By Fact~\ref{Sylow_fusion}, $N_M(S^\o)$ controls $M$-fusion in
$C_S(S^\o)$, so all involutions in $E$ are conjugate in $N_H(S^\o)$.
Hence $E \leqs \Omega_1(S^\o)$ and $\pr_2(S) \geqs 3$, in
contradiction with Lemma~\ref{nonconjugacy}.
\end{verification}

\begin{Claim}\label{stremb_Gmin}
$G$ is a minimal connected simple group.
\end{Claim}

\begin{verification}
Suppose towards a contradiction that $G$ has a proper definable
 non-solvable connected subgroup.
Let $K$ be a minimal proper definable non-solvable connected
subgroup of $G$ and let $\-K=K/O(K)$.  By Fact~\ref{511}, $\-K$ is
a central product of quasi-simple algebraic groups over algebraically closed
fields of characteristic not 2 and of one definable connected abelian
group.  By minimality of $K$, $\-K$ must actually be one quasi-simple
algebraic group.  Now $\-K$ must be isomorphic to either $\SL_2(F)$ or
$\PSL_2(F)$ for some algebraically closed field $F$ of characteristic
not 2, also by minimality of $K$ (Fact~\ref{alggrpmin}).

A 2-generated core of $K$ is a subgroup of a 2-generated core of $G$.
So the connected component of a 2-generated core of $K$ is also solvable.
By Lemma~\ref{lifting}, the connected component of a 2-generated core
of $\-K$ is also solvable.
By Lemma~\ref{notpsl2}, $\-K \not\simeq \PSL_2(F)$, so
$\-K \simeq \SL_2(F)$.

Now $\-K$ has a central involution $\bar z$ in the connected component
of a Sylow 2-subgroup.
By Fact \ref{Cquotient}, $C_K(z)O(K)/O(K) = C_{K/O(K)}(\bar z)$
 for some involution $z\in S^\o$.
By Claim~\ref{stremb_CconM} applied to $E$, $C_G^\o(z) \leqs B$,
in contradiction with \clref{stremb_Msolv}.
\end{verification}

Now suppose for the moment that $S$ is connected.  Then $S$ is abelian
and $C_G(i) \leq M$ for every $i\in M$ by \clref{stremb_CM}.
Hence $M$ is strongly embedded by Fact~\ref{strembdef}.
% ; and the proof of Theorem~\ref{stremb} is complete.
Since $\pr_2(S) = n_2(S) \geq 3$ too, $B$ is a Borel subgroup
by Lemma~\ref{twogencore_max_pr}.  This means that we can dedicate the
remainder of the argument to showing that $S$ is connected.

\begin{Claim}\label{stremb_NBM}
$N_G(B) = M$  %%% (so $[N_G(B):B] < \infty$).
\end{Claim}

\begin{verification}
We observe that $S^\o$ is now a Sylow 2-subgroup of $B$ by Fact~\ref{Sylow_con}.
Since $\pr_2(S) \geq 2$, $N_G(S^\o) \leq N_G(\Omega_1(S^\o)) \leq M$.
By Fact~\ref{Sylow_frattini},
\[ N_G(B) = B N_{N_G(B)}(S^\o) \leq M \qedhere \]
\end{verification}

\begin{Claim}\label{stremb_IBBg}
$I(B \cap B^g) = \emptyset$ for any $g\not\in N_G(B)$.
\end{Claim}

\begin{verification}
Suppose towards a contradiction that there is an $i\in I(B \cap B^g)$.
We may assume that $i \in I(S)$ after conjugation. % by Fact~\ref{Sylow}.
Since $B$ is solvable by \clref{stremb_Msolv}, the Sylow
2-subgroups of $B$ are connected by Fact~\ref{Sylow_con}.
As $i\in B$ and $B$ has odd type, $S^\circ \leq C^\o_G(i)$.
Since $i\in M^g$, $C^\o_G(i) \leq M^g$ by \clref{stremb_CconM}.
Since $\pr_2(S) \geq 2$, $S \leq N_G(\Omega_1(S^\o)) \leq M^g$.
Since $S$ is now a Sylow 2-subgroup of $M^g$,
$M^g = \Gamma_{S,2}(G) = M$.
\end{verification}

\begin{Claim}\label{stremb_BGgen}
$\bigcup B^G$ is generic in $G$.
\end{Claim}

For this, we employ the following fact from \cite{CJ01}
 and a general lemma.

\begin{fact}[Lemma 3.3 of \cite{CJ01}]\label{stremb_BGgen_L}
Let $G$ be a connected group of finite Morley rank.  Let $B$ be a definable
subgroup of finite index in its normalizer.  Suppose there is a definable
subset $Q$ of $B$, not generic in $B$, such that $B \cap B^g \subseteq Q$
whenever $g\notin N_G(B)$.  Then $\bigcup B^G$ is generic in $G$.
\end{fact}

\begin{lemma}\label{solvable_2perp_union}
Let $H$ be a connected solvable group of finite Morley rank and odd type.
Let $\mathcal F$ be a uniformly definable family of $2^\perp$-subgroups
of $H$.  Then there is a definably characteristic definable
$2^\perp$-subgroup $Q$ of $H$ containing $\bigcup \mathcal F$.
\end{lemma}

Here definably characteristic means invariant under definable automorphisms.

\begin{proof}%[Proof of Lemma \ref{solvable_2perp_union}]
Lemma 3.2 of \cite{CJ01} says that the quotient $\bar{H} := H/O(H)$
 is divisible abelian, since $H$ is connected solvable of odd type.
So $\bar{F} \normal \bar{H}$ for any $F\in \mathcal F$.
By Fact~\ref{nodeftorsion}, $\bar{F}$ is a $2^\perp$-subgroup of $\bar{H}$
for any $F\in \mathcal F$.  Since $O(\bar{H}) = 1$ and $\bar{H}$ is abelian,
$\bar{F}$ is finite for any $F\in \mathcal F$.
Since the family $\{ \bar{F} : F\in \mathcal F \}$ is uniformly definable,
there is a bound on $|\bar{F}|$ by Axiom D of \cite[p.~57]{BN}.
So $m = \lcm \{ |\bar{F}| : F\in \mathcal F \} < \infty$ is odd.
Since $\bar{H}$ is abelian, $\bar{Q} := \{ h\in \bar{H} : h^m = 1 \}$
is a characteristic $2^\perp$-subgroup of $\bar{H}$ containing
$\bar{F}$ for all $F\in \mathcal F$.  So the pullback $Q$ of $\bar{Q}$ is
a suitable definably characteristic $2^\perp$-subgroup of $H$.
% by Fact~\ref{nodeftorsion}.
\end{proof}

\begin{verification}[Verification of \clref{stremb_BGgen}]
By \clref{stremb_IBBg} and Lemma~\ref{solvable_2perp_union},
there is a definably characteristic definable $2^\perp$-subgroup
$Q$ of $B$ which contains $B \cap B^g$ for any $g\notin N_G(B)$.
Since $B$ has non-trivial Sylow 2-subgroup, we have $Q < B$.
$B$ has finite index in its normalizer by \clref{stremb_NBM}.
Now $\bigcup B^G$ is generic in $G$ by Fact~\ref{stremb_BGgen_L}.
\end{verification}

We observe that conclusions 5, 6, and 7 follow from the previous three claims.

Consider a pair $B_1,B_2$ of definable subgroups of $G$.  We say a definable subgroup $H$ of $G$ is
{\it $(B_1,B_2)$-bi-invariant} if $H$ is $(A_1,A_2)$-invariant for some four-groups $A_1 \leq B_1$
and $A_2 \leq B_2$.   We may simply say $H$ is {\it bi-invariant} when the choice of $B_1$ and $B_2$
are clear from the context.  Similarly, we say that a collection of definable subgroups $\mathcal H$ is
{\it simultaneously bi-invariant} if all $H\in \mathcal H$ are $(A_1,A_2)$-invariant for the same
choice of $A_1$ and $A_2$.

\begin{Claim}\label{stremb_bi_inv}
Any connected definable $(B,B^g)$-bi-invariant subgroup $K$ of $G$ is contained in $B \cap B^g$;
and hence is a $2^\perp$-group when $g\notin N_G(B)$ by \clref{stremb_IBBg}.
\end{Claim}

\begin{verification}
By Fact~\ref{Kgrp_gen} and \clref{stremb_CconM}
\[ K = \langle C_K^\o(a) : a\in A_1 \rangle \leq M \]
and similarly $K \leq M^g$.
% \[ K = \langle C_K^\o(a) : a\in A_2 \rangle \leq M^g
%  \qedhere \]
\end{verification}

We claim that $S$ is connected.  Suppose towards a contradiction that $S$ is disconnected.
We fix an $i\in S-S^\o$ with $i^2 \in S^\o$. We also define
$$ X := \{ x \in iB : \textrm{$x \in (\gen{i} B)^g$
                      for some $g\notin N_G(B)$} \} $$

\begin{Claim}\label{stremb_j_exists}
There is a $j\in X$ with $j^2 = 1$.
\end{Claim}

For this, we employ the following fact from \cite{CJ01}.

\begin{fact}[Lemma 3.5 of \cite{CJ01}]\label{stremb_Xgen_L}
Let $G$ be a connected group of finite Morley rank.  Let $B$ be a proper
definable subgroup of finite index in its normalizer such that $\bigcup B^G$
is generic in $G$.  Suppose that $z\in N_G(B)-B$ has order $n>1$ modulo $B$,
and let $\gen{z} B$ be the union
 $B \cup z B \cup z^2 B \cup \cdots \cup z^{n-1} B$.
Then the following subset $X$ of $z B$ is generic in $z B$.
$$ X := \{ x \in z B :
 \textrm{$x \in (\gen{z} B)^g$ for some $g\notin N_G(B)$} \} $$
\end{fact}

\begin{verification}[Verification of \clref{stremb_j_exists}]
$X$ is generic in $iB$ by Fact~\ref{stremb_Xgen_L} and \clref{stremb_BGgen}.
So there is some $x\in X$.  Then $x\in iB \cap (\gen{i} B)^g$ for some
$g\notin N_G(B)$ and $x^2\in B \cap B^g$.  So $K := \{1,x\}(B \cap B^g)$
is a definable group, and $B \cap B^g \normal K$. % since $[K:B \cap B^g] = 2$.
%%% Should I say that $K4 is either $\gen{i}B \cap (\gen{i}B)^g$ or $\gen{i}B \cap B^g$?
By Fact~\ref{nodeftorsion}, there is a non-trivial 2-element
$j\in x(B \cap B^g) \leq X$.  Now $j^2=1$ since $j^2\in B \cap B^g$
and $I(B \cap B^g) = \emptyset$ by \clref{stremb_IBBg}.
\end{verification}

For the next portion of our argument, we fix an involution $j\in X$
and some $g\notin N_G(B)$ with $j\in (\gen{i} B)^g$.
% We consider the intersection of $B$ and $B^g$.

\begin{Claim}\label{stremb_Cconj}
% There are four-groups $A_1 \leq E$ and $A_2 \leq E^g$ such that
$C^\o_G(j)$ is non-trivial and is $(B,B^g)$-bi-invariant.
\end{Claim}

\begin{verification}
Since $G$ is non-abelian, $C^\o_G(j)$ is non-trivial \cite[Ex.~13 p.~79]{BN}
% by Fact~\ref{centinv}. % by Exercise 13 on page 79 of \cite{BN} (from \cite{Ne90a}).
Since $j\in iB \leq M$, $j\in S^b$ for some $b\in M$ by conjugacy. % Fact~\ref{Sylow}.
By \clref{stremb_E1}, there is a four-group $A_1\leq E^b$ centralizing
$j$.  The existence of a suitable $A_2 \leq B^g$ follows similarly.
% by a simillar argument with $M^G$ and $E^g$.
\end{verification}

Now consider a maximal proper definable connected $(B,B^g)$-bi-invariant
subgroup $H$ of $G$.  $H$ is non-trivial by \clref{stremb_Cconj}.
Let $C$ be a Carter subgroup of $H$ (which exists by Fact~\ref{Carter_exists}-1).

 From this point forward, we will have no more need of the assumption
that $S$ is disconnected, or the involution $j$.  Instead, we proceed
by general arguments involving the groups $B$, $B^g$, $H$, and $C$.

\begin{Claim}\label{stremb_C_inv}
$C$ and $H$ are simultaneously $(B,B^g)$-bi-invariant.
\end{Claim}

\begin{verification}
Let $A$ denote one of the two groups with respect to which $H$ is bi-invariant,
i.e.~either $A_1$ or $A_2$.  By Fact~\ref{Carter_frattini}-3, $H A = H N_{HA}(C)$.
Since $I(H) = \emptyset$ by \clref{stremb_bi_inv}, $A$ is a Sylow 2-subgroup of
$H A$ and $N_{H A}(C)/C \cong H A / H$ is a four-group.
By Fact~\ref{Sylow_lift}, $N_{H A}(C)$ contains a Sylow 2-subgroup $A_0$ of $HA$.
So $H$ and $C$ are clearly $A_0$-invariant and $A_0$ lives in $B$ or $B^g$, as appropriate.
\end{verification}

\begin{Claim}\label{stremb_Carter}
$N^\o_G(C) = C$; and hence $C$ is a Carter subgroup of any Borel
subgroup containing it.
\end{Claim}

\begin{verification}
The two groups $H$ and $N^\o_G(C)$ are simultaneously bi-invariant by \clref{stremb_C_inv},
so $\gen{ H, N^\o_G(C) } \leq B \cap B^g$ is proper and bi-invariant.
Hence $N^\o_G(C) \leq H$ by maximality.  So $N^\o_G(C) = C$; and
$C$ is a Carter subgroup of any Borel subgroup which contains it
by Fact~\ref{Carter_cleanup}-2.
\end{verification}

The following general lemma now shows that $C$ contains a conjugate of $S^\o$.

\begin{lemma}\label{Carter_Sylow_cor}
Let $B$ be a connected solvable group of finite Morley rank and odd type,
and let $C$ be a Carter subgroup of $B$.  Then $C$ contains a conjugate
of the Sylow 2-subgroup of $B$.
\end{lemma}

\begin{proof}
Let $S$ be a Sylow 2-subgroup of $B$.  By Fact~\ref{Carter_Sylow}-4,
$N_B(S)$ contains a Carter subgroup $C_1$ of $B$.
Since $C_1$ is connected by Fact~\ref{Carter_con}-2,
$C_1 \leqs N^\o_B(S) \leqs C_B(S)$ by Lemma 6.16 of \cite{BN}.
% By Fact~\ref{Sylow_con}, $S$ is connected.
So $S \leqs N_B(C_1) = C_1$.  By Fact~\ref{Carter_conj}-3,
$C_1$ is conjugate to $C$.
\end{proof}

By Lemma \ref{Carter_Sylow_cor} below, $C$ contains a conjugate of $S^\o$;
in contradiction to \clref{stremb_bi_inv}.  Thus $S$ is connected and all of
our claims follow. \qed

% We observe that Claims \ref{stremb_NBM}, \ref{stremb_IBBg}, and
% \ref{stremb_BGgen}, are free of contradictory assumptions, and we may
% extend the Strong Embedding theorem by these additional claims:
% \begin{enumerate}
% \item[5]  $N_G(B) = M$.
% \item[6]  $I(B \cap B^g) = \emptyset$ for any $g\notin M$.
% \item[7]  $\bigcup B^G$ is generic in $G$.
% \end{enumerate}

\subsection*{Acknowledgments}

The authors thank Ay\c{s}e Berkman for stimulating discussions and
Tuna \Altinel, Gregory Cherlin and Eric Jaligot for helpful comments
and corrections.

% \nocite{Al94,ACCN98,Bo98,Carter89,kleidman-liebeck} % wagner
\small

\def\cprime{$'$}

\end{document}